\documentclass[11pt]{article}
\usepackage{mathrsfs,amsmath,amssymb,amsfonts}
\numberwithin{equation}{section}
\begin{document}
\title{Pointwise weighted approximation of functions with inner singularities by Bernstein operators}
\author{Wen-ming Lu$^1$, Lin Zhang$^2$\thanks{\textsl{E-mail address}:\texttt{lu\_wenming@163.com}(W.-Lu).}}
\date{\small \it $^1$School of Science, Hangzhou Dianzi University,
Hangzhou, 310018 P.R. China\\$^2$Department of Mathematics, Zhejiang University, Hangzhou, 310027
P.R. China} \maketitle
\mbox{}\hrule\mbox{}\\[0.5cm]
\textbf{Abstract}\\[-0.2cm]^^L
We consider the pointwise weighted approximation by Bernstein
operators with inner singularities. The related weight functions are
weights $\bar{w}(x)=|x-\xi|^\alpha(0<\xi<1,\ \alpha>0).$ In this
paper we give direct and inverse results of this type of Bernstein
polynomials.\\~\\
\textbf{Keywords:} Pointwise weighted approximation; Bernstein
operators; inner singularities\\[0.5cm]
\mbox{}\hrule\mbox{}
\section{Introduction}
The set of all continuous functions, defined on the interval $I$, is
denoted by $C(I)$. For any $f\in C([0,1])$, the corresponding
\emph{Bernstein operators} are defined as follows:
$$B_n(f,x):=\sum_{k=0}^nf(\frac{k}{n})p_{n,k}(x),$$
where
$$p_{n,k}(x)={n \choose k}x^k(1-x)^{n-k}, \ k=0,1,2,\ldots,n, \ x\in[0,1].$$
Approximation properties of Bernstein operators have been studied
very well. Berens and Lorentz showed in \cite{Berens} that
\begin{eqnarray*}
B_n(f,x)-f(x)=O(({\frac{1}{\sqrt{n}}}\delta_n(x))^{\alpha_1})
\Longleftrightarrow \omega^1(f,t)=O(t^{\alpha_1}),
\end{eqnarray*}
where $0< \alpha_1 <1,\
\delta_n(x)=\varphi(x)+{\frac{1}{\sqrt{n}}},\
\varphi(x)=\sqrt{x(1-x)}.$\\~\\
It is well known that approximation of functions with singularities
by polynomial is of special value in both theories and applications.
As an important type of polynomial approximation, approximation of
functions by Bernstein operators is an important topic in both
approximation theory and computational theory, which plays an
important role in neural networks, fitting date, curves, and
surfaces. Some work has been done by \cite{Della}. Throughout the
paper, $C$ denotes a positive constant independent of $n$ and $x$,
which may be different in different cases.\\~\\
Let ${\bar{w}}(x)=|x-\xi|^\alpha,\ 0<\xi<1,\ \alpha>0$ and
$C_{\bar{w}}:= \{{f \in C([0,1] \setminus \{\xi\})
:\lim\limits_{x\longrightarrow\xi}(\bar{w}f)(x)=0 }\}$. The
\emph{norm} in $C_{\bar{w}}$ is defined by
$\|f\|_{C_{\bar{w}}}:=\|{\bar{w}}f\|=\sup\limits_{0\leqslant
x\leqslant 1}|({\bar{w}}f)(x)|$. Define
$$W_{\bar{w},\lambda}^{2}:= \{f \in
C_{\bar{w}}:f' \in A.C.((0,1)),\
\|{\bar{w}}\varphi^{2\lambda}f''\|<\infty\}.$$ For $f \in
C_{\bar{w}}$, the \emph{weighted modulus of smoothness} is defined
by
\begin{eqnarray*}
\omega_{\varphi}^{2}(f,t)_{\bar{w}}:=\sup_{0<h\leqslant
t}\{\|{\bar{w}}\triangle_{h\varphi}^{2}f\|_{[16h^2,1-16h^2]}+\|{\bar{w}}{\overrightarrow{\triangle}_{h}^{2}}f\|_{[0,16h^2]}+\|{\bar{w}}{\overleftarrow{\triangle}_{h}^{2}}f\|_{[1-16h^2,1]}\},
\end{eqnarray*}
where
\begin{eqnarray*}
\Delta_{h\varphi}^{2}f(x)&=&f(x+h\varphi(x))-2f(x)+f(x-h\varphi(x)),\\
{\overrightarrow{\triangle}_{h}^{2}}f(x)&=&f(x+2h)-2f(x+h)+f(x),\\
{\overleftarrow{\triangle}_{h}^{2}}f(x)&=&f(x-2h)-2f(x-h)+f(x),
\end{eqnarray*}
and $\varphi(x)=\sqrt{x(1-x)},\
\delta_n(x)=\varphi(x)+{\frac{1}{\sqrt{n}}}.$\\~\\
Let
\begin{eqnarray*}
\psi(x)=\left\{
\begin{array}{lrr}
10x^3-15x^4+6x^5, &&0< x <1, \\
0,   &&x \leqslant0,  \\
1,  &&x \geqslant 1.
             \end{array}
\right.
\end{eqnarray*}
Obviously, $\psi$ is non-decreasing on the real axis, $\psi\in
C^2((-\infty,+\infty)),\ \psi^{(i)}(0)=0,$ $i=0,1,2.\
\psi^{(i)}(1)=0,\ i=1,2$ and $\psi(1)=1.$ Further, let
\begin{eqnarray*}
x_{1}=\frac{[n\xi-2\sqrt{n}]}{n},\ x_{2}=\frac{[n\xi-\sqrt{n}]}{n},\ x_{3}=\frac{[n\xi+\sqrt{n}]}{n},\ x_{4}=\frac{[n\xi+2\sqrt{n}]}{n},
\end{eqnarray*}
and
\begin{eqnarray*}
{\bar{\psi}}_{1}(x)=\psi(\frac{x-x_{1}}{x_{2}-x_{1}}),\ {\bar{\psi}}_{2}(x)=\psi(\frac{x-x_{3}}{x_{4}-x_{3}}).
\end{eqnarray*}
Consider
\begin{eqnarray*}
P(x):={\frac{x-x_{4}}{x_{1}-x_{4}}}f(x_{1}) + {\frac{x_{1}-x}{x_{1}-x_{4}}}f(x_{4}),
\end{eqnarray*}
the linear function joining the points $(x_{1},f(x_{1}))$ and
$(x_{4},f(x_{4}))$. And let
\begin{eqnarray*}
{\bar{F}}_{n}(f,x):={\bar{F}}_{n}(x)=f(x)(1-{\bar{\psi}}_{1}(x)+{\bar{\psi}}_{2}(x))+{\bar{\psi}}_{1}(x)(1-{\bar{\psi}}_{2}(x))P(x).
\end{eqnarray*}
From the above definitions it follows that
\begin{eqnarray*}
{\bar{F}}_{n}(f,x)=\left\{\begin{array}{lr}
f(x),          &       x\in [0,x_1]\cup [x_4,1],   \\
f(x)(1-{\bar{\psi}}_{1}(x))+{\bar{\psi}}_{1}(x)
P(x),      &
x\in [x_1,x_2],  \\
P(x),          &       x\in [x_2,x_3],  \\
P(x)(1-{\bar{\psi}}_{2}(x))+{\bar{\psi}}_{2}(x)f(x), & x\in
[x_3,x_4].
            \end{array}
\right.
\end{eqnarray*}
Evidently, ${\bar{F}}_{n}$ is a positive linear operator which
depends on the functions values $f(k/n),\ 0 \leqslant k/n \leqslant
x_2$ or $x_3 \leqslant k/n \leqslant 1,$ it reproduces linear
functions, and ${\bar{F}}_{n} \in C^2([0,1])$ provided $f\in
W_{\bar{w},\lambda}^{2}.$ Now for every $f\in C_{\bar{w}}$ define
the \emph{Bernstein type operator}
\begin{eqnarray}
{\bar{B}}_{n}(f,x)&:=&B_{n}({\bar{F}}_{n}(f),x)  \nonumber\\
&=&\sum_{k/n \in [0,x_1]\cup [x_4,1]}p_{n,k}(x)f({\frac kn}) + \sum_{x_2<k/n<x_3}p_{n,k}(x)P({\frac kn})\nonumber\\
&&+\sum_{x_1<k/n<x_2}p_{n,k}(x)\{f({\frac kn})(1-{\bar{\psi}}_{1}({\frac kn}))+{\bar{\psi}}_{1}({\frac kn})P({\frac kn})\}\nonumber\\
&&+\sum_{x_3<k/n<x_4}p_{n,k}(x)\{P({\frac
kn})(1-{\bar{\psi}}_{2}({\frac kn}))+{\bar{\psi}}_{2}({\frac
kn})f({\frac kn})\} \label{s1}
\end{eqnarray}
Obviously, ${\bar{B}}_{n}$ is a positive linear operator,
${\bar{B}}_{n}(f)$ is a polynomial of degree at most $n,$ it
preserves linear functions, and depends only on the function values
$f(k/n),\ k/n \in [0,x_2] \cup [x_3,1].$
Now we state our main results as follows:\\~\\
\textbf{Theorem 1.} {\it If \ $\alpha>0,$ for any $f \in
C_{\bar{w}},$ we have $\|\bar{w}{\bar{B}}''_n(f)\|\leqslant
Cn^2\|{\bar{w}}f\|$.} \\~\\
\textbf{Theorem 2.} {\it For any $\alpha
>0,\ 0 \leqslant \lambda \leqslant 1,$ we have
\begin{eqnarray*}
|{\bar{w}(x)}\varphi^{2\lambda}(x){\bar{B}}''_n(f,x)|\leqslant
\left\{
\begin{array}{lrr}
Cn{\{max\{n^{1-\lambda},\varphi^{2(\lambda-1)}\}\}}\|{\bar{w}}f\|,    &&f\in C_{\bar{w}},    \\
C\|{\bar{w}}\varphi^{2\lambda}f''\|, &&f\in
W_{{\bar{w}},\lambda}^{2}.
              \end{array}
\right.
\end{eqnarray*}}
\textbf{Theorem 3.} {\it For $f\in C_{\bar{w}},\ 0<\xi<1,\ \alpha>0,\ \alpha_0 \in(0,2),$ we have
\begin{eqnarray*}
{\bar{w}(x)}|f(x)-{\bar{B}_{n}(f,x)}|=O((n^{-{\frac 12}}\varphi^{-\lambda}(x)\delta_n(x))^{\alpha_0}) \Longleftrightarrow \omega_{\varphi^\lambda}^2(f,t)_{\bar{w}}=O(t^{\alpha_0}).
\end{eqnarray*}}
\section{Lemmas}
\textbf{Lemma 1.}(\cite{Zhou}) {\it For any non-negative real $u$ and
$v$, we have
\begin{eqnarray}
\sum_{k=1}^{n-1}({\frac kn})^{-u}(1-{\frac
kn})^{-v}p_{n,k}(x)\leqslant Cx^{-u}(1-x)^{-v}. \label{s2}
\end{eqnarray}}
\textbf{Lemma 2.}(\cite{Della}) {\it For any $\alpha >0,\ f \in
C_{\bar{w}},$ we have
\begin{eqnarray}
\|{\bar{w}}{\bar{B}_{n}(f)}\| \leqslant C\|{\bar{w}}f\|. \label{s3}
\end{eqnarray}}
\textbf{Lemma 3.}(\cite{Wang}) {\it  If \
$\varphi(x)=\sqrt{x(1-x)},\ 0 \leqslant \lambda \leqslant 1,\ 0
\leqslant \beta \leqslant 1,$ then
\begin{eqnarray}
\int_{-{\frac {h\varphi^\lambda(x)}{2}}}^{\frac
{h\varphi^\lambda(x)}{2}} \cdots \int_{-{\frac
{h\varphi^\lambda(x)}{2}}}^{\frac
{h\varphi^\lambda(x)}{2}}\varphi^{-r\beta}(x+\sum_{k=1}^ru_k)du_1
\cdots du_r \leqslant Ch^r\varphi^{r(\lambda-\beta)}(x).\label{s4}
\end{eqnarray}}
\textbf{Lemma 4.}(\cite{Della}) {\it If $\gamma \in R,$ then
\begin{eqnarray}
\sum_{k=0}^np_{n,k}(x)|k-nx|^\gamma \leqslant Cn^{\frac
\gamma2}\varphi^\gamma(x).\label{s5}
\end{eqnarray}}
\textbf{Lemma 5.}  {\it Let
$A_n(x):={\bar{w}(x)}\sum\limits_{|k-n\xi|\leqslant
\sqrt{n}}p_{n,k}(x)$. Then $A_n(x)\leqslant Cn^{-\alpha/2}$ for
$0<\xi <1$ and $\alpha>0$.} \\~\\
\textbf{Proof.} If $|x- \xi|\leqslant {\frac {3}{\sqrt{n}} }$, then
the statement is trivial. Hence assume $0 \leqslant x \leqslant
\xi-{\frac {3}{\sqrt{n}} }$ (the case $\xi+{\frac {3}{\sqrt{n}} }
\leqslant x \leqslant 1$ can be treated similarly). Then for a fixed
$x$ the maximum of $p_{n,k}(x)$ is attained for
$k=k_n:=[n\xi-\sqrt{n}]$. By using Stirling's formula, we get
\begin{eqnarray*}
p_{n,k_n}(x)&\leqslant& C{\frac {({\frac
{n}{e}})^n\sqrt{n}x^{k_n}(1-x)^{n-k_n}}{({\frac
{k_n}{e}})^{k_n}\sqrt{k_n}({\frac
{n-k_n}{e}})^{n-k_n}\sqrt{n-k_n}}}\nonumber\\
&\leqslant& {\frac {C}{\sqrt{n}}}({\frac {nx}{k_n}})^{k_n}({\frac
{n(1-x)}{n-k_n}})^{n-k_n}\nonumber\\
&=&{\frac {C}{\sqrt{n}}}(1-{\frac {k_n-nx}{k_n}})^{k_n}(1+{\frac
{k_n-nx}{n-k_n}})^{n-k_n}.
\end{eqnarray*}
Now from the inequalities
$$k_n-nx=[n\xi-\sqrt{n}]-nx>n(\xi-x)-\sqrt{n}-1\geqslant {\frac 12}n(\xi-x),$$
and
$$1-u\leqslant e^{-u-{\frac 12}u^2},\ 1+u\leqslant e^u,\ u\geqslant 0,$$
it follows that the second inequality is valid. To prove the first
one we consider the function $\lambda(u)=e^{-u-{\frac 12}u^2}+u-1.$
Here $\lambda(0)=0,\ \lambda^\prime(u)=-(1+u)e^{-u-{\frac
12}u^2}+1,\ \lambda^\prime(0)=0,\
\lambda^{\prime\prime}(u)=u(u+2)e^{-u-{\frac 12}u^2}\geqslant 0,$
whence $\lambda(u)\geqslant 0$ for $u\geqslant 0$. Hence
\begin{eqnarray*}
p_{n,k_n}(x) &\leqslant&{\frac {C}{\sqrt{n}}}exp\{k_n[-{\frac
{k_n-nx}{k_n}}-{\frac 12}({\frac {k_n-nx}{k_n}})^2] +
k_n-nx\}\nonumber\\
&=&{\frac {C}{\sqrt{n}}}exp\{-\frac {({k_n-nx})^2}{2k_n}\}\leqslant
e^{-Cn(\xi-x)^2}.
\end{eqnarray*}
Thus $A_n(x)\leqslant C(\xi-x)^\alpha e^{-Cn(\xi-x)^2}$. An easy
calculation shows that here the maximum is attained when
$\xi-x={\frac {C}{\sqrt{n}}}$ and the lemma follows. $\Box$\\~\\
\textbf{Lemma 6.} {\it For $0<\xi <1,\ \alpha,\ \beta>0$, we have
\begin{eqnarray}
{\bar{w}(x)}\sum\limits_{|k-n\xi|\leqslant \sqrt{n}}|k-nx|^\beta
p_{n,k}(x)\leqslant Cn^{(\beta-\alpha)/2}\varphi^\beta(x).\label{s6}
\end{eqnarray}}
\textbf{Proof.} By (\ref{s5}) and the lemma 5, we have
\begin{eqnarray*}
{\bar{w}(x)}^{\frac
{1}{2n}}({\bar{w}(x)\sum\limits_{|k-n\xi|\leqslant
\sqrt{n}}p_{n,k}(x)})^{\frac
{2n-1}{2n}}(\sum\limits_{|k-n\xi|\leqslant \sqrt{n}}|k-nx|^{2n\beta}
p_{n,k}(x))^{\frac {1}{2n}}\leqslant
Cn^{(\beta-\alpha)/2}\varphi^\beta(x). \Box
\end{eqnarray*}
\textbf{Lemma 7.} {\it For any $\alpha >0,\ f\in
W_{\bar{w},\lambda}^{2},$ we have
\begin{eqnarray}
{\bar{w}(x)}|f(x)-P(f,x)|_{[x_1,x_4]} \leqslant C({\frac
{\delta_n(x)}{\sqrt{n}\varphi^\lambda(x)}})^2\|{\bar{w}}\varphi^{2\lambda}f''\|.\label{s7}
\end{eqnarray}}
\textbf{Proof.} If $x \in [x_1,x_4],$ for any $f\in
W_{\bar{w},\lambda}^{2},$ we have
\begin{eqnarray*}
f(x_1)=f(x)+f'(x)(x_1-x)+\int_{x_1}^x (t-x_1)f''(t)dt,\\
f(x_4)=f(x)+f'(x)(x_4-x)+\int_{x_4}^x (t-x_4)f''(t)dt,\\
{\delta_n(x)} \sim {\frac {1}{\sqrt{n}}},\ n=1,2,\cdots.
\end{eqnarray*}
So
\begin{eqnarray*}
{\bar{w}(x)}|f(x)-P(f,x)| &\leqslant& {\bar{w}(x)}|{\frac
{x-x_4}{x_1-x_4}}|\int_{x_1}^x
|(t-x_1)f''(t)|dt\\
&&+{\bar{w}(x)}|{\frac
{x_1-x}{x_1-x_4}}|\int_{x_4}^x |(t-x_4)f''(t)|dt\\
&:=&I_1+I_2.
\end{eqnarray*}
Whence $t$ is between $x_1$ and $x$. We have
$\frac{|t-x_1|}{{\bar{w}}(t)}\leqslant
\frac{|x-x_1|}{{\bar{w}}(x)},$ then
\begin{eqnarray*}
I_1 &\leqslant&
C\|{\bar{w}}\varphi^{2\lambda}f''\||(x-x_1)(x-x_4)|\int_{x_1}^x\varphi^{-2\lambda}(t)dt\\
&\leqslant&
C({\frac{\delta_n(x)}{\sqrt{n}\varphi^\lambda(x)}})^2\|{\bar{w}}\varphi^{2\lambda}f''\|.
\end{eqnarray*}
Analogously, we have
\begin{eqnarray*}
I_2 \leqslant
C({\frac{\delta_n(x)}{\sqrt{n}\varphi^\lambda(x)}})^2\|{\bar{w}}\varphi^{2\lambda}f''\|.
\end{eqnarray*}
Now the lemma follows from combining these results together. $\Box$
\section{Proof of Theorem}
\subsection{Proof of Theorem 1}
If $f \in C_{\bar{w}},$ when $x\in [{\frac 1n},1-{\frac 1n}],$ by
[2], we have
\begin{eqnarray*}
|{\bar{w}(x)}{\bar{B}''_{n}(f,x)}| & \leqslant &
n\varphi^{-2}(x){\bar{w}(x)}|{\bar{B}_{n}(f,x)}|\\
&+&{\bar{w}(x)}\varphi^{-4}(x)\sum_{k=0}^np_{n,k}(x)|k-nx||{\bar{F}}_{n}({\frac
kn})|\\
& + &
{\bar{w}(x)}\varphi^{-4}(x)\sum_{k=0}^n(k-nx)^2|{\bar{F}}_{n}({\frac
kn})|p_{n,k}(x)\\  &:=& A_1+A_2+A_3.
\end{eqnarray*}
By (\ref{s3}), we have
\begin{eqnarray*}
A_1(x)=n\varphi^{-2}(x){\bar{w}(x)}|{\bar{B}_{n}(f,x)}| \leqslant
Cn^2\|{\bar{w}}f\|.
\end{eqnarray*}
and
\begin{eqnarray*}
A_2&=&{\bar{w}(x)}\varphi^{-4}(x)[\sum_{k/n\in
A}|k-nx||{\bar{F}}_{n}({\frac kn})|p_{n,k}(x)+\sum_{x_2 \leqslant
k/n\leqslant x_3}|k-nx||P({\frac
kn})|p_{n,k}(x)]\\&:=&\sigma_1+\sigma_2.
\end{eqnarray*}
thereof $A:=[0,x_2]\cup [x_3,1].$ If ${\frac kn}\in A,$ when ${\frac
{\bar{w}(x)}{\bar{w}(\frac {k}{n})}}\leqslant C(1+n^{-{\frac
{\alpha}{2}}}|k-nx|^\alpha),$ we have $|k-n\xi|\geqslant {\frac
{\sqrt{n}}{2}}$, by (\ref{s5}), then
\begin{eqnarray*}
\sigma_1 & \leqslant &
C\|{\bar{w}}f\|\varphi^{-4}(x)\sum_{k=0}^np_{n,k}(x)|k-nx|[1+n^{-{\frac
{\alpha}{2}}}|k-nx|^\alpha]\\
&=&
C\|{\bar{w}}f\|\varphi^{-4}(x)\sum_{k=0}^np_{n,k}(x)|k-nx|+Cn^{-{\frac
{\alpha}{2}}}\|{\bar{w}}f\|\varphi^{-4}(x)\sum_{k=0}^np_{n,k}(x)|k-nx|^{1+\alpha}\\
& \leqslant & Cn^{\frac 12}\varphi^{-3}(x)\|{\bar{w}}f\|+Cn^{\frac
12}\varphi^{-3+\alpha}(x)\|{\bar{w}}f\| \\
& \leqslant & Cn^2\|{\bar{w}}f\|.
\end{eqnarray*}
For $\sigma_2,$ $P$ is a linear function. We note $|P({\frac
kn})|\leqslant max(|P(x_1)|,|P(x_4)|):=P(a).$ If $x\in [x_1,x_4],$
we have ${\bar{w}(x)}\leqslant {\bar{w}(a)}.$ So, if $x\in
[x_1,x_4],$ by (\ref{s5}), then
\begin{eqnarray*}
\sigma_2\leqslant
C{\bar{w}(a)}P(a)\varphi^{-4}(x)\sum_{k=0}^np_{n,k}(x)|k-nx|\leqslant
Cn^{2}\|{\bar{w}}f\|.
\end{eqnarray*}
If $x\notin [x_1,x_4],$ then ${\bar{w}(a)}>n^{-{\frac
{\alpha}{2}}},$ by (\ref{s6}), we have
\begin{eqnarray*}
\sigma_2 &\leqslant& C{\bar{w}(x)}\varphi^{-4}(x) \sum_{x_2
\leqslant k/n\leqslant x_3}|P(a)(k-nx)|p_{n,k}(x)\\
&\leqslant& Cn^{\frac
{\alpha}{2}}\|{\bar{w}}f\|{\bar{w}(x)}\varphi^{-4}(x)\sum_{x_2
\leqslant k/n\leqslant x_3}|k-nx|p_{n,k}(x)\\
&\leqslant& Cn^2\|{\bar{w}}f\|.
\end{eqnarray*}
So, $A_2 \leqslant Cn^2\|{\bar{w}}f\|$. Similarly, $A_3 \leqslant
Cn^2\|{\bar{w}}f\|$. It follows from combining the above
inequalities that the inequality is proved.\\~\\
When $x\in
[0,{\frac 1n}]$ (The same as $x\in [1-{\frac 1n},1]$), by
\cite{Lorentz}, then
\begin{eqnarray*}
{\bar{B}''_{n}(f,x)}=n(n-1)\sum_{k=0}^{n-2}\overrightarrow{\Delta}_{\frac
1{n}}^{2}{\bar{F}}_{n}{(\frac kn)}p_{n-2,k}(x).
\end{eqnarray*}
We have
\begin{eqnarray*}
|{\bar{w}(x)}{\bar{B}''_{n}(f,x)}| &\leqslant&
Cn^2{\bar{w}(x)}\sum_{k=0}^{n-2}|\overrightarrow{\Delta}_{\frac
1{n}}^{2}{\bar{F}}_{n}{(\frac kn)}|p_{n-2,k}(x)\\
&=& Cn^2{\bar{w}(x)}[\sum_{k/n\in
A}p_{n-2,k}(x)|\overrightarrow{\Delta}_{\frac
1{n}}^{2}{\bar{F}}_{n}{(\frac kn)}|+\sum_{x_2 \leqslant k/n\leqslant
x_3}p_{n-2,k}(x)|\overrightarrow{\Delta}_{\frac 1{n}}^{2}P({\frac
kn})|].
\end{eqnarray*}
We can deal with it in accordance with the former proofs, and prove
it immediately, then the theorem is done. $\Box$
\subsection{Proof of Theorem 2}
(1) We prove the first inequality of Theorem 2.\\~\\
\textit{Case 1.} If $0\leqslant \varphi(x)\leqslant {\frac
{1}{\sqrt{n}}}$, by Theorem 1, we have
\begin{eqnarray*}
|{\bar{w}(x)}\varphi^{2\lambda}(x){\bar{B}}''_n(f,x)|\leqslant
Cn^{-\lambda}|{\bar{w}(x)}{\bar{B}}''_n(f,x)|\leqslant
Cn^{2-\lambda}\|{\bar{w}}f\|.
\end{eqnarray*}
\textit{Case 2.} If $\varphi(x)> {\frac {1}{\sqrt{n}}}$, by
\cite{Totik},\ we have
\begin{eqnarray*}
{\bar{B}}''_n(f,x)=B''_n({\bar{F}_{n}},x) & = &
(\varphi^{2}(x))^{-2}\sum_{i=0}^{2}Q_{i}(x,n)n^i\sum_{k=0}^{n}(x-{\frac
kn})^{i}{\bar{F}}_{n}({\frac kn})p_{n,k}(x),\\
Q_{i}(x,n) &=& (nx(1-x))^{[(2r-i)/2]},\\
(\varphi^{2}(x))^{-2}Q_{i}(x,n)n^{i} &\leqslant&
C(n/\varphi^{2}(x))^{1+i/2}.
\end{eqnarray*}
So
\begin{eqnarray*}
&&|{\bar{w}}(x)\varphi^{2\lambda}(x){\bar{B}}''_n(f,x)|\nonumber\\
&\leqslant& C{\bar{w}(x)}\varphi^{2\lambda}(x)\sum_{i=0}^{2}({\frac
{n}{\varphi^2(x)}})^{1+i/2}\sum_{k=0}^{n}|(x-{\frac
kn})^{i}{\bar{F}}_{n}({\frac kn})|p_{n,k}(x)\nonumber\\
&=& C{\bar{w}(x)}\varphi^{2\lambda}(x)\sum_{i=0}^{2}({\frac
{n}{\varphi^2(x)}})^{1+i/2}\sum_{k/n\in A}|(x-{\frac
kn})^{i}{\bar{F}}_{n}({\frac kn})|p_{n,k}(x)\nonumber\\
&&+ C{\bar{w}(x)}\varphi^{2\lambda}(x)\sum_{i=0}^{2}({\frac
{n}{\varphi^2(x)}})^{1+i/2}\sum_{x_2 \leqslant k/n \leqslant
x_3}|{(x-{\frac kn})^{i}}P({\frac
kn})|p_{n,k}(x)\nonumber\\
&:=&\sigma_1+ \sigma_2,
\end{eqnarray*}
where $A:=[0,x_2]\cup [x_3,1]$. Working as in the proof of Theorem
1, we can get $\sigma_1\leqslant Cn^{2-\lambda}\|{\bar{w}}f\|,$
$\sigma_2\leqslant Cn^{2-\lambda}\|{\bar{w}}f\|.$ By bringing these
facts together, we can immediately get the first inequality of
Theorem 2.\\~\\
(2) If $f\in W_{\bar{w},\lambda}^{2},$ by
${\bar{B}}_{n}(f,x)=B_{n}({\bar{F}}_{n}(f),x),$ then
\begin{eqnarray*}
|{\bar{w}(x)}\varphi^{2\lambda}(x){\bar{B}}''_n(f,x)| & \leqslant &
n^2{\bar{w}(x)}\varphi^{2\lambda}(x)\sum_{k=0}^{n-2}|\overrightarrow{\Delta}_{\frac
1{n}}^{2}{\bar{F}}_{n}{(\frac kn)}|p_{n-2,k}(x)\\
&=&
n^2{\bar{w}(x)}\varphi^{2\lambda}(x)\sum_{k=1}^{n-3}|\overrightarrow{\Delta}_{\frac
1{n}}^{2}{\bar{F}}_{n}{(\frac kn)}|p_{n-2,k}(x)\\
&+&
n^2{\bar{w}(x)}\varphi^{2\lambda}(x)|\overrightarrow{\Delta}_{\frac
1{n}}^{2}{\bar{F}}_{n}{(0)}|p_{n-2,0}(x)\\
&+&
n^2{\bar{w}(x)}\varphi^{2\lambda}(x)|\overrightarrow{\Delta}_{\frac
1{n}}^{2}{\bar{F}}_{n}{(\frac {n-2}{n})}|p_{n-2,n-2}(x) \\
&:=& I_1+I_2+I_3.
\end{eqnarray*}
By \cite{Totik}, if $0<k<n-2,$ we have
\begin{eqnarray}
|\overrightarrow{\Delta}_{\frac 1{n}}^{2}{\bar{F}}_{n}({\frac kn})
|\leqslant Cn^{-1}\int_{0}^{\frac {2}{n}}|{\bar{F}}''_{n}({\frac
kn}+u)|du,\label{s8}
\end{eqnarray}
If $k=0,$ we have
\begin{eqnarray}
|\overrightarrow{\Delta}_{\frac 1{n}}^{2}{\bar{F}}_{n}(0)| \leqslant
C\int_{0}^{\frac {2}{n}}u|{\bar{F}}''_{n}(u)|du,\label{s9}
\end{eqnarray}
Similarly
\begin{eqnarray}
|\overrightarrow{\Delta}_{\frac 1n}^2{\bar{F}}_{n}({\frac
{n-2}{n}})| &\leqslant& Cn^{-1}\int_{1-{\frac
{2}{n}}}^{1}(1-u)|{\bar{F}}''_{n}(u)|du.\label{s10}
\end{eqnarray}
By (\ref{s8}), then
\begin{eqnarray}
I_1 &\leqslant&
Cn{\bar{w}(x)}\varphi^{2\lambda}(x)\sum_{k=1}^{n-3}\int_{0}^{\frac
{2}{n}}|{\bar{F}}''_{n}({\frac kn}+u)|dup_{n-2,k}(x) \nonumber\\
&=& Cn{\bar{w}(x)}\varphi^{2\lambda}(x)\sum_{k/n\in
A}\int_{0}^{\frac {2}{n}}|{\bar{F}}''_{n}({\frac
kn}+u)|dup_{n-2,k}(x) \nonumber\\
&&+ Cn{\bar{w}(x)}\varphi^{2\lambda}(x)\sum_{x_2 \leqslant k/n
\leqslant x_3}\int_{0}^{\frac {2}{n}}|P''({\frac
kn}+u)|dup_{n-2,k}(x) \nonumber\\
&:=& T_1+T_2,\label{s11}
\end{eqnarray}
where $A:=[0,x_2]\cup [x_3,1],$ $P$ is a linear function. If ${\frac
kn}\in A,$ when ${\frac {\bar{w}(x)}{\bar{w}(\frac
{k}{n})}}\leqslant C(1+n^{-{\frac {\alpha}{2}}}|k-nx|^\alpha),$ we
have $|k-n\xi|\geqslant {\frac {\sqrt{n}}{2}}$, by (\ref{s2}),
(\ref{s5}) and the Theorem 2, then
\begin{eqnarray*}
T_1 &\leqslant&
C{\bar{w}(x)}\varphi^{2\lambda}(x)\|{\bar{w}}\varphi^{2\lambda}{\bar{F}''_{n}}\|\sum_{k/n\in
A}p_{n-2,k}(x){\bar{w}}^{-1}({\frac kn})\varphi^{-2\lambda}({\frac kn})\\
&\leqslant&
C\varphi^{2\lambda}(x)\|{\bar{w}}\varphi^{2\lambda}{\bar{F}''_{n}}\|\sum_{k=0}^{n-2}p_{n-2,k}(x)(1+n^{-{\frac
{\alpha}{2}}}|k-nx|^\alpha)\varphi^{-2\lambda}({\frac kn})\\
&\leqslant& C\|{\bar{w}}\varphi^{2\lambda}{\bar{F}''_{n}}\|\\
&\leqslant& C\|{\bar{w}}\varphi^{2\lambda}f''\|.
\end{eqnarray*}
Working as the Theorem 1, we can get
\begin{eqnarray*}
T_2 \leqslant C\|{\bar{w}}\varphi^{2\lambda}f''\|.
\end{eqnarray*}
So, we can get
\begin{eqnarray*}
I_1 \leqslant C\|{\bar{w}}\varphi^{2\lambda}f''\|.
\end{eqnarray*}
By (\ref{s9}) and the Theorem 2, we have
\begin{eqnarray}
I_2 &\leqslant&
Cn^2{\bar{w}(x)}\varphi^{2\lambda}(x)(1-x)^{n-2}\int_{0}^{\frac
{2}{n}}u|{\bar{F}}''_{n}(u)|du  \nonumber\\
&\leqslant&
Cn^2{\bar{w}(x)}\varphi^{2\lambda}(x)(1-x)^{n-2}\|{\bar{w}}\varphi^{2\lambda}{\bar{F}''_{n}}\|\int_{0}^{\frac
{2}{n}}u{\bar{w}}^{-1}(u)\varphi^{-2\lambda}(u)du  \nonumber\\
&\leqslant& C\|{\bar{w}}\varphi^{2\lambda}{\bar{F}''_{n}}\|  \nonumber\\
&\leqslant& C\|{\bar{w}}\varphi^{2\lambda}f''\|.  \label{s12}
\end{eqnarray}
Similarly,
\begin{eqnarray}
I_3 \leqslant C\|{\bar{w}}\varphi^{2\lambda}f''\|.\label{s13}
\end{eqnarray}
By bringing (\ref{s11}), (\ref{s12}) and (\ref{s13}) together, we
can get the second inequality
of Theorem 2. $\Box$\\~\\
\textbf{Corollary 1.} {\it If $\alpha>0$ and $\lambda=0$, we have
\begin{eqnarray*}
|{\bar{w}(x)}{\bar{B}}''_n(f,x)|\leqslant \left\{
\begin{array}{lrr}
Cn^2\|{\bar{w}}f\|,    &&f\in C_{\bar{w}},    \\
C\|{\bar{w}}f''\|, &&f\in W_{{\bar{w}}}^2.
              \end{array}
\right.
\end{eqnarray*}}
\textbf{Corollary 2.} {\it If $\alpha>0$ and $\lambda=1$, we have
\begin{eqnarray*}
|{\bar{w}(x)}\varphi^2(x){\bar{B}}''_n(f,x)|\leqslant \left\{
\begin{array}{lrr}
Cn\|{\bar{w}}f\|,    &&f\in C_{\bar{w}},   \\
C\|{\bar{w}}\varphi^{2}f''\|, &&f\in W_{{\bar{w}}}^2.
\end{array} \right.
\end{eqnarray*}}
\subsection{Proof of Theorem 3}
\subsubsection{The direct theorem}
We know
\begin{eqnarray*}
{\bar{F}}_n(t)={\bar{F}}_n(x)+{\bar{F}}'_n(t)(t-x)+\int_x^t (t-u){\bar{F}}''_n(u)du, \\
B_n(t-x,x)=0.
\end{eqnarray*}
According to the definition of $W_{\bar {w},\lambda}^{2},$ \  for
any $g \in W_{\bar {w},\lambda}^{2},$ we have
${\bar{B}}_{n}(g,x)=B_{n}({\bar{G}}_{n}(g),x)$.\\~\\
(1) We first
estimate ${\bar{w}(x)}|{\bar{G}}_{n}(x)-B_{n}({\bar{G}}_{n},x)|$
under the condition of $x\in [{\frac 1n},1-{\frac 1n}],$ then
$\varphi^2(x)<{\frac 1n},\ \delta_n(x) \sim {\frac {1}{\sqrt{n}}}$,
and
\begin{eqnarray*}
{\bar{w}(x)}|{\bar{G}}_{n}(x)-B_{n}({\bar{G}}_{n},x)|={\bar{w}(x)}|B_n(R_2({\bar{G}}_n,t,x),x)|
\end{eqnarray*}
thereof $R_2({\bar{G}}_n,t,x)=\int_x^t (t-u){\bar{G}}''_n(u)du.$
\\
It follows from $\frac{|t-u|}{{\bar{w}}(u)}\leqslant
\frac{|t-x|}{{\bar{w}}(x)},$ $u$ between $t$ and $x$, we have
\begin{eqnarray*}
{\bar{w}(x)}|{\bar{G}}_{n}(x)-B_{n}({\bar{G}}_{n},x)| & \leqslant &
C\|{\bar{w}}\varphi^{2\lambda}{\bar{G}}''_n\|{\bar{w}(x)}B_n(\int_x^t{\frac
{|t-u|}{{\bar{w}(u)}\varphi^{2\lambda}(u)}du,x})\\
& \leqslant &
C\|{\bar{w}}\varphi^{2\lambda}{\bar{G}}''_n\|{\bar{w}}(x)(B_n(\int_x^t{\frac
{|t-u|}{\varphi^{4\lambda}(u)}}|du,x))^{\frac 12}(B_n(\int_x^t{\frac
{|t-u|}{{\bar{w}^2(u)}}}du,x))^{\frac 12}
\end{eqnarray*}
also
\begin{eqnarray}
\int_x^t{\frac {|t-u|}{\varphi^{4\lambda}(u)}}du \leqslant C{\frac
{(t-x)^2}{\varphi^{4\lambda}(x)}},\ \int_x^t{\frac
{|t-u|}{{\bar{w}^2(u)}}}du \leqslant {\frac
{(t-x)^2}{{\bar{w}^2(x)}}}. \label{s14}
\end{eqnarray}
By (\ref{s5}) and (\ref{s14}), we have
\begin{eqnarray*}
{\bar{w}(x)}|{\bar{G}}_{n}(x)-B_{n}({\bar{G}}_{n},x)| &\leqslant&
C\|{\bar{w}}\varphi^{2\lambda}{\bar{G}}''_n\|\varphi^{-2\lambda}(x)B_n
((t-x)^2,x)\\
&\leqslant& Cn^{-1}{\frac
{\varphi^{2}(x)}{\varphi^{2\lambda}(x)}}\|{\bar{w}}\varphi^{2\lambda}{\bar{G}}''_n\|\\
&\leqslant& Cn^{-1}{\frac
{\delta_n^{2}(x)}{\varphi^{2\lambda}(x)}}\|{\bar{w}}\varphi^{2\lambda}{\bar{G}}''_n\|\\
&=& C(\frac
{\delta_n(x)}{\sqrt{n}\varphi^\lambda(x)})^2\|{\bar{w}}\varphi^{2\lambda}{\bar{G}}''_n\|
\end{eqnarray*}
(2) We estimate
${\bar{w}(x)}|{\bar{G}}_{n}(x)-B_{n}({\bar{G}}_{n},x)|$ under the
condition of $x\in [0,{\frac 1n})$ (The same as $x\in (1-{\frac
1n},1]$), $\varphi(x) \sim \delta_n(x),$  now
\begin{eqnarray*}
{\bar{w}(x)}|{\bar{G}}_{n}(x)-B_{n}({\bar{G}}_{n},x)| &\leqslant&
C{\bar{w}(x)}\sum_{k=1}^{n-1}p_{n,k}(x)\int_x^{\frac kn}|({\frac
kn}-u){\bar{G}}''_n(u)|du\\
&+& C{\bar{w}(x)}p_{n,0}(x)\int_0^xu|{\bar{G}}''_n(u)|du\\
&+&
C{\bar{w}(x)}p_{n,n}(x)\int_x^1|(1-u){\bar{G}}''_n(u)|du\\
&:=& I_1+I_2+I_3.
\end{eqnarray*}
If $u$ between ${\frac kn}$ and $x,$ we have
\begin{eqnarray}
\frac{|{\frac kn}-u|}{{\bar{w}}^2(u)}\leqslant \frac{|{\frac
kn}-x|}{{\bar{w}}^2(x)},\
\frac{|k/n-u|}{\varphi^{4\lambda}(u)}\leqslant
\frac{|k/n-x|}{\varphi^{4\lambda}(x)}.\label{s15}
\end{eqnarray}
By (\ref{s5}) and (\ref{s15}), then
\begin{eqnarray}
I_1 &\leqslant&
C\|{\bar{w}}\varphi^{2\lambda}{\bar{G}}''_n\|{\bar{w}(x)}\sum_{k=1}^{n-1}p_{n,k}\int_x^{\frac
kn}{\frac {|k/n-u|}{{\bar{w}(u)}\varphi^{2\lambda}(u)}}du  \nonumber\\
&\leqslant&
C\|{\bar{w}}\varphi^{2\lambda}{\bar{G}}''_n\|{\bar{w}(x)}\sum_{k=1}^{n-1}p_{n,k}(\int_x^{\frac
kn}{\frac {|k/n-u|}{\bar{w}^2(u)}}du)^{\frac 12}({\int_x^{\frac
kn}{\frac
{|k/n-u|}{\varphi^{4\lambda}(u)}}du})^{\frac 12}  \nonumber\\
&\leqslant&Cn^{-2}\|{\bar{w}}\varphi^{2\lambda}{\bar{G}}''_n\|\varphi^{-2\lambda}(x)\sum_{k=0}^{n-1}p_{n,k}(x)(k-nx)^2  \nonumber\\
&\leqslant& Cn^{-1}{\frac
{\varphi^{2}(x)}{\varphi^{2\lambda}(x)}}\|{\bar{w}}\varphi^{2\lambda}{\bar{G}}''_n\| \nonumber\\
&\leqslant& Cn^{-1}{\frac
{\delta_n^{2}(x)}{\varphi^{2\lambda}(x)}}\|{\bar{w}}\varphi^{2\lambda}{\bar{G}}''_n\|  \nonumber\\
&=& C(\frac
{\delta_n(x)}{\sqrt{n}\varphi^\lambda(x)})^2\|{\bar{w}}\varphi^{2\lambda}{\bar{G}}''_n\|.
\label{s16}
\end{eqnarray}
For $I_2,$ when $u$ between ${\frac kn}$ and $x,$ we let $k=0,$ then
$\frac{u}{{\bar{w}}(u)}\leqslant \frac{x}{{\bar{w}}(x)},$ and
\begin{eqnarray}
I_2 &\leqslant&
C\|{\bar{w}}\varphi^{2\lambda}{\bar{G}}''_n\|{\bar{w}(x)}p_{n,0}(x)\int_0^xu{\bar{w}^{-1}(u)}\varphi^{-2\lambda}(u)du\nonumber\\
&\leqslant& C(nx)(1-x)^{n-1} \cdot n^{-1}{\frac
{\varphi^{2}(x)}{\varphi^{2\lambda}(x)}}\|{\bar{w}}\varphi^{2\lambda}{\bar{G}}''_n\|\nonumber\\
&\leqslant& C(\frac
{\delta_n(x)}{\sqrt{n}\varphi^\lambda(x)})^2\|{\bar{w}}\varphi^{2\lambda}{\bar{G}}''_n\|.\label{s17}
\end{eqnarray}
Similarly, we have
\begin{eqnarray}
I_3 \leqslant C(\frac
{\delta_n(x)}{\sqrt{n}\varphi^\lambda(x)})^2\|{\bar{w}}\varphi^{2\lambda}{\bar{G}}''_n\|.\label{s18}
\end{eqnarray}
By bringing (\ref{s16}), (\ref{s17}) and (\ref{s18}), we get the
result. Above all, we have
\begin{eqnarray*}
{\bar{w}(x)}|{\bar{G}}_{n}(x)-B_{n}({\bar{G}}_{n},x)| \leqslant
C(\frac
{\delta_n(x)}{\sqrt{n}\varphi^\lambda(x)})^2\|{\bar{w}}\varphi^{2\lambda}{\bar{G}}''_n\|.
\end{eqnarray*}
By (\ref{s7}) and the second inequality of Theorem 2, when $g \in
W_{\bar{w},\lambda}^{2},$ then
\begin{eqnarray}
{\bar{w}(x)}|g(x)-{\bar{B}_{n}(g,x)}| &\leqslant&
{\bar{w}(x)}|g(x)-{\bar{G}}_{n}(g,x)| +
{\bar{w}(x)}|{\bar{G}}_{n}(g,x)-{\bar{B}_{n}(g,x)}|\nonumber\\
&\leqslant& {\bar{w}(x)}|g(x)-P(g,x)|_{[x_1,x_4]} + C(\frac
{\delta_n(x)}{\sqrt{n}\varphi^\lambda(x)})^2\|{\bar{w}}\varphi^{2\lambda}{\bar{G}}''_n\|\nonumber\\
&\leqslant& C(\frac
{\delta_n(x)}{\sqrt{n}\varphi^\lambda(x)})^2\|{\bar{w}}\varphi^{2\lambda}g''\|.\label{s19}
\end{eqnarray}
For $f \in C_{\bar{w}},$ we choose proper $g \in W_{\bar
{w},\lambda}^{2},$ by (\ref{s3}) and (\ref{s19}), then
\begin{eqnarray*}
{\bar{w}(x)}|f(x)-{\bar{B}_{n}(f,x)}| &\leqslant&
{\bar{w}(x)}|f(x)-g(x)| + {\bar{w}(x)}|{\bar{B}_{n}(f-g,x)}| +
{\bar{w}(x)}|g(x)-{\bar{B}_{n}(g,x)}|\\
&\leqslant& C(\|{\bar{w}}(f-g)\|+(\frac
{\delta_n(x)}{\sqrt{n}\varphi^\lambda(x)})^2\|{\bar{w}}\varphi^{2\lambda}g''\|)\\
&\leqslant& C\omega_{\varphi^\lambda}^2(f,\frac
{\delta_n(x)}{\sqrt{n}\varphi^\lambda(x)})_{\bar{w}}. \Box
\end{eqnarray*}
\subsubsection{The inverse theorem}
We define the weighted main-part modulus for $D=R_+$ by
\begin{eqnarray*}
\Omega_{\varphi^\lambda}^2(C,f,t)_{\bar{w}} = \sup_{0<h \leqslant
t}\|{\bar{w}}\Delta_{{h\varphi}^\lambda}^2f\|_{[Ch^\ast,\infty]},\\
\Omega_{\varphi^\lambda}^2(1,f,t)_{\bar{w}} =
\Omega_{\varphi^\lambda}^2(f,t)_{\bar{w}}.
\end{eqnarray*}
where $C>2^{1/\beta(0)-1},\ \beta(0)>0,$ and $h^\ast$ is given by
\begin{eqnarray*}
h^\ast= \left\{
\begin{array}{lrr}
(Ar)^{1/1-\beta(0)}h^{1/1-\beta(0)},  && 0 \leqslant \beta(0) <1,
    \\
0,   && \beta(0) \geqslant 1.
              \end{array}
\right.
\end{eqnarray*}
The main-part $K$-functional is given by
\begin{eqnarray*}
H_{\varphi^\lambda}^2(f,t^2)_{\bar{w}}=\sup_{0<h \leqslant
t}\inf_g\{\|{\bar{w}}(f-g)\|_{[Ch^\ast,\infty]}+t^2\|{\bar{w}}\varphi^{2\lambda}g''\|_{[Ch^\ast,\infty]},\
g' \in A.C.((Ch^\ast,\infty))\}.
\end{eqnarray*}
By \cite{Totik}, we have
\begin{eqnarray}
C^{-1}\Omega_{\varphi^\lambda}^2(f,t)_{\bar{w}} \leqslant
\omega_{\varphi^\lambda}^{2}(f,t)_{\bar{w}} \leqslant
C\int_0^t{\frac {\Omega_{\varphi^\lambda}^2(f,\tau)_{\bar{w}}}{\tau}}d\tau,\label{s20} \\
C^{-1}H_{\varphi^\lambda}^{2}(f,t^2)_{\bar{w}} \leqslant
\Omega_{\varphi^\lambda}^2(f,t)_{\bar{w}} \leqslant
CH_{\varphi^\lambda}^{2}(f,t^2)_{\bar{w}}.\label{s21}
\end{eqnarray}
\textbf{Proof.} Let $\delta>0,$ by (\ref{s21}), we choose proper $g$
so that
\begin{eqnarray*}
\|{\bar{w}}(f-g)\| \leqslant
C\Omega_{\varphi^\lambda}^2(f,t)_{\bar{w}},\
\|{\bar{w}}\varphi^{2\lambda}g''\| \leqslant
C\delta^{-2}\Omega_{\varphi^\lambda}^2(f,t)_{\bar{w}}.
\end{eqnarray*}
then
\begin{eqnarray}
|{\bar{w}}(x)\Delta_{h\varphi^\lambda}^2f(x)| &\leqslant&
|{\bar{w}}(x)\Delta_{h\varphi^\lambda}^2(f(x)-{\bar{B}_{n}(f,x)})|+|{\bar{w}}(x)\Delta_{h\varphi^\lambda}^2\bar{B}_{n}(f-g,x)|\nonumber\\
&&+ |{\bar{w}}(x)\Delta_{h\varphi^\lambda}^2{\bar{B}_{n}(g,x)}|\nonumber\\
&\leqslant& \sum_{j=0}^2C_2^j(n^{-\frac
12}\delta_n(x+(1-j)h\varphi^\lambda(x)))^{\alpha_0}\nonumber\\
&&+ \int_{-{\frac {h\varphi^\lambda(x)}{2}}}^{\frac
{h\varphi^\lambda(x)}{2}} \int_{-{\frac
{h\varphi^\lambda(x)}{2}}}^{\frac
{h\varphi^\lambda(x)}{2}}{\bar{w}}(x){\bar{B}''_{n}(f-g,x+\sum_{k=1}^2u_k)}du_1
du_2\nonumber\\
&&+ \int_{-{\frac {h\varphi^\lambda(x)}{2}}}^{\frac
{h\varphi^\lambda(x)}{2}} \int_{-{\frac
{h\varphi^\lambda(x)}{2}}}^{\frac
{h\varphi^\lambda(x)}{2}}{\bar{w}}(x){\bar{B}''_{n}(g,x+\sum_{k=1}^2u_k)}du_1
du_2\nonumber\\
&:=& J_1+J_2+J_3.\label{s22}
\end{eqnarray}
Obviously
\begin{eqnarray}
J_1 \leqslant C(n^{-\frac 12}\delta_n(x))^{\alpha_0}.\label{s23}
\end{eqnarray}
By Theorem 1, we have
\begin{eqnarray}
J_2 &\leqslant& Cn^2\|{\bar{w}}(f-g)\|\int_{-{\frac
{h\varphi^\lambda(x)}{2}}}^{\frac {h\varphi^\lambda(x)}{2}}
\int_{-{\frac {h\varphi^\lambda(x)}{2}}}^{\frac
{h\varphi^\lambda(x)}{2}}du_1du_2\nonumber\\
&\leqslant& Cn^2h^2\varphi^{2\lambda}(x)\|{\bar{w}}(f-g)\|\nonumber\\
&\leqslant&
Cn^2h^2\varphi^{2\lambda}(x)\Omega_{\varphi^\lambda}^2(f,\delta)_{\bar{w}}.\label{s24}
\end{eqnarray}
By the second inequality of Corollary 2 and (\ref{s4}),\ we have
\begin{eqnarray}
J_2 &\leqslant& Cn\|{\bar{w}}(f-g)\|\int_{-{\frac
{h\varphi^\lambda(x)}{2}}}^{\frac {h\varphi^\lambda(x)}{2}}
\int_{-{\frac {h\varphi^\lambda(x)}{2}}}^{\frac
{h\varphi^\lambda(x)}{2}}\varphi^{-2}(x+\sum_{k=1}^2u_k)du_1du_2\nonumber\\
&\leqslant& Cnh^2\varphi^{2(\lambda-1)}(x)\|{\bar{w}}(f-g)\|\nonumber\\
&\leqslant&
Cnh^2\varphi^{2(\lambda-1)}(x)\Omega_{\varphi^\lambda}^2(f,\delta)_{\bar{w}}.\label{s25}
\end{eqnarray}
By the second inequality of Theorem 2 and (\ref{s4}),\ we have
\begin{eqnarray}
J_3 &\leqslant&
C\|{\bar{w}}\varphi^{2\lambda}g''\|{\bar{w}(x)}\int_{-{\frac
{h\varphi^\lambda(x)}{2}}}^{\frac {h\varphi^\lambda(x)}{2}}
\int_{-{\frac {h\varphi^\lambda(x)}{2}}}^{\frac
{h\varphi^\lambda(x)}{2}}{\bar{w}^{-1}(x+\sum_{k=1}^2u_k)}\varphi^{-2\lambda}(x+\sum_{k=1}^2u_k)du_1du_2\nonumber\\
&\leqslant& Ch^2\|{\bar{w}}\varphi^{2\lambda}g''\|\nonumber\\
&\leqslant&
Ch^2\delta^{-2}\Omega_{\varphi^\lambda}^2(f,\delta)_{\bar{w}}.\label{s26}
\end{eqnarray}
Now, by (\ref{s23}), (\ref{s24}), (\ref{s25}) and (\ref{s26}), we
get
\begin{eqnarray*}
|{\bar{w}}(x)\Delta_{h\varphi^\lambda}^2f(x)| &\leqslant&
C\{(n^{-\frac 12}\delta_n(x))^{\alpha_0} + h^2(n^{-\frac
12}\delta_n(x))^{-2}\Omega_{\varphi^\lambda}^2(f,\delta)_{\bar{w}} +
h^2\delta^{-2}\Omega_{\varphi^\lambda}^2(f,\delta)_{\bar{w}}\}.
\end{eqnarray*}
When $n \geqslant 2,$ we have
\begin{eqnarray*}
n^{-\frac 12}\delta_n(x) < (n-1)^{-\frac 12}\delta_{n-1}(x)
\leqslant \sqrt{2}n^{-\frac 12}\delta_n(x),
\end{eqnarray*}
Choosing proper $x, n \in N,$ so that
\begin{eqnarray*}
n^{-\frac 12}\delta_n(x) \leqslant \delta < (n-1)^{-\frac
12}\delta_{n-1}(x),
\end{eqnarray*}
Therefore
\begin{eqnarray*}
|{\bar{w}}(x)\Delta_{h\varphi^\lambda}^2f(x)| \leqslant
C\{\delta^{\alpha_0} +
h^2\delta^{-2}\Omega_{\varphi^\lambda}^2(f,\delta)_{\bar{w}}\}.
\end{eqnarray*}
By Borens-Lorentz lemma, we get
\begin{eqnarray}
\Omega_{\varphi^\lambda}^2(f,t)_{\bar{w}} \leqslant
Ct^{\alpha_0}.\label{s27}
\end{eqnarray}
So, by (\ref{s27}), we get
\begin{eqnarray*}
\omega_{\varphi^\lambda}^{2}(f,t)_{\bar{w}} \leqslant
C\int_0^t{\frac
{\Omega_{\varphi^\lambda}^2(f,\tau)_{\bar{w}}}{\tau}}d\tau =
C\int_0^t\tau^{\alpha_0-1}d\tau = Ct^{\alpha_0}. \Box
\end{eqnarray*}
\subsection*{Acknowledgement}
~\indent The authors would like to thank the anonymous referees
whose comments have been implemented in the final version of the
manuscript.

\end{document}